\documentclass[12pt]{article}
\usepackage{amsmath}
\usepackage{amssymb}
\begin{document}
\begin{center}

{\bf Dense set of large periodic points}\\
\vspace{.2in}
S.Kanmani\\
Material Science Division, IGCAR,\\ Kalpakkam, India.\\
\vspace{.1in}
V.Kannan\\
Department of Mathematics and Statistics,\\
University of Hyderabad, Hyderabad.\\
India.\\

\end{center}

{\bf Abstract- } We say that a dynamical system $(X,f)$ has a
dense set of large periodic points  (abbreviated as $DLP$)
if for every positive integer $n$, the set of all periodic points
of period $ \geq n $, is dense in $X$. Here we provide two different
proofs of the fact that every system chaotic in the sense of Devaney,
has $DLP$.

 {\bf Keywords - } Devaney's chaos, periodic orbits, Topological transitivity

\subsection{Introduction}

A dynamical system is a pair $ (X,f) $, where $X$ is a Hausdorff space
and $f$ is a continuous self-map of $X$. The set of all positive integers
is denoted by $N$. An element $x  \in  X $ is said to be periodic if
$f^{n}(x)=x $ for some $n \in N $; the smallest such $n$ is called the
period of $x$. Here $f^n$ denotes the $n$-fold composition of $ f $. The
set of all periodic points is denoted $ P $. For each $n \in N, P_{n} $
denotes the set of all periodic  points of period $ \leq n$. 

{\bf Definition}  We say that the system has $DLP$ if $P \setminus P_{n}$ is dense in
X for every $n$ in $N$.
According to Devaney [1] $ (X,f) $ is said to be chaotic if it has all
the three properties below.

\begin{enumerate}

\item
The set of all periodic points $( P ) $ is dense in X.
\item
$ f $ is topologically transitive and

\item
$f $ has sensitive dependence on initial conditions.

\end{enumerate}
But it is known [ 2] that $ 3 $ folows from $ 1 $ and $ 2 $.
We say that $ f $ is topologically transitive if for any two non-empty 
open sets $ V $ and $ W $ in $X$, there exist $n$ in $N$ such that
$ f^{n}(V) \cap W $ is not empty. Roughly speaking, the requirement is
that from any region $ V $ to any region $ W $, some point will go at
some time.  

	The main result in this note is that every Devaney chaotic system
possesses $ DLP $. In the field of study like  Monte carlo methods one
uses a discrete dynamical system as a pseudo-random number generator.
There one of the desirable features of the dynamical system is to have
orbits whose period  is as high as possible. Our result gives a
criteria which ensures the abundance of large periodic points. 

	It is known [3, 4] that any continuous function on an interval
of the real line, if topologically transitive,  has  infinitely
many periodic orbits, and that the periods of these orbits forms an unbounded subset
of the positive integers. Our main result allows to draw a stronger conclusion
from the same hypotheses. This is possible because it is known on intervals topological 
transitivity  implies that the set of all periodic points is dense [5].

	Our arguements prove that the property of
dense set of periodic points and  dense set of large periodic points
are equivalent among  the following two classes of dynamical systems :

\begin{enumerate}

\item
 Topologically transitive dynamical systems that are infinite.
 
\item

Dynamical systems in which  the set of all periodic points has an
empty interior.

\end{enumerate}

However in general,  dense $P$ and $DPL$ are not equivalent.
Because, a dense $P$ does not necessarily imply $DLP$. Some simple examples are:

\begin{enumerate}

\item

 Cyclic permutations on finite sets. Here $ P = $ Whole set $ X $ and  $P \setminus  P_{n}$ is empty 
if  $n$ exceeds the cardinality of $ X $.

\item
Identity map on any $ X $. Here $ P= X $ and $P \setminus P_{n}$ is empty for $ n \geq 2 $

\item
Reflection maps on the real line $R$. Here $ P = R $ and $P \backslash P_{n}$ is empty if
$n \geq 3 $. (e.g. $f(x)= 1-x $)

\end{enumerate}

We shall later note that every continuous map on $R$  that has dense set of periodic points but
without dense large periodic points, must have a subsystem that is topologically conjugate
to a reflection or identity.

\subsection{The two proofs}

The main theorem is:

{\bf Theorem-1 } Every transitive map with dense set of periodic points has $DLP$, if the underlying space is 
infinite.

{\bf Proof: }

	Our first proof draws this as an easy corollary of the main result of Touhey [6]. We say 
that a collection $ V_{1},V_{2} \ldots  , V_{n} $ of nonempty open sets, shares a periodic orbit if
there is a periodic point $x$ whose orbit meets each $V_{j} , 1\leq j \leq n $. Here, the orbit of $x$
is the set of all elements of the form $f^n(x)$ where $n$ is in $N$. It follows from [6] that every 
transitive map with dense set of periodic points  has the  property that every finite collection of nonempty
 open sets shares a
periodic
orbit.

Since $X$ is infinite and $f$ is transitive, it is easily seen that $X$ has no isolated points. Next, if
$V$ is any open set in $X$, then $V$ is infinite (as there are no isolated points). Let $ n \in N $. Since
$X$ is a Hausdorff space $n+1$ points  in $V$ can be mutually separated by disjoint open sets. That is, we
can find a family $ V_{1},V_{2} \ldots  , V_{n+1} $ of pairwise disjoint non-empty open sets inside $V$.

By applying the main result of [6], we get a periodic point $x$ whose orbit meets each $V_{i}$. Because
these are pairwise disjoint, it is immediate that this orbit has atleast $n+1$ elements. If $y$ is an element
in this orbit that is also in $V_{1}$, then $y$ is a periodic point in $V$ whose period is strictly
greater than $n$. This proves that $ P \setminus P_{n} $ is dense in $X$.

	The second proof uses the following facts:

\begin{enumerate}

\item
Each $P_{n}$ is closed.

\item
Each $P_{n}$ is invariant.

\item
In a topologally transitive system, every invariant set is either dense or nowhere dense.

\item
When a nowhere dense set is removed from a dense set, what remains is again a dense set.

\end{enumerate}

We consider two cases:

Case-1:  $ P_{n} $ is not equal to the whole space $X$ for any $n$. 
Then by $ (1) $ and $ (3) $ each $P_{n}$ is nowhere dense. By (4) $P \setminus P_{n} $ is dense. Thus the
system has the $DPL$ property.

Case-2: $ P_{n}=X $ for some $n$. 
Since the dynamical  system is transitive,
first we observe that if $V$ and $W$ are two nonempty
open sets, then  $ f^{m}(V) \cap W $ is
 non-empty for some $ m < n $. It follows from the proposition $1$ below that $X$ has atmost
$n$ elements. Since $X$ is finite, transitivity demands that the
dynamical system is  a cyclic permutation i.e.
$f^{n} $ is the identity map.
\subsection{Speciality of $n$-cycles:}

	Topological transitivity requires that from every nonempty open set $V$ it
is possible to reach every nonempty open set $W$ at some time or other. Can we put
an upperbound to the time of reaching ? That is, can we demand that there is a positive
integer $k$ such that $f^n(V) \cap W$ is nonempty for some  $n < k $, whatever $V$ and
$W$ be ? That this is not possible, except in the trivial cases, is the contention of the next
proposition.

 {\bf Proposition-1 } Let $(X,f)$ be a dynamical system and let $k$ be positive integer. Suppose that
whenever $V$ and $W$ are non-empty open sets in $X$, there is a positive integer $ n < 
k$
such that $f^n(V) \cap W $ is non-empty. Then $X$ has atmost $k$ elements.

Proof: Let if possible, $ x_{0}, x_{1}, \ldots x_{k} $ be $k+1$ distinct elements
 of $X$. Choose pairwise disjoint
neighbourhoods $ V_{0}, V_{1}, \ldots V_{k} $ of these points. This is possible if $X$ is a metric space, or more generally
$X$
is a Hausdorff topological space. Then choose a neighbourhood $W_{0}$ of $  x_{0}$ and
a positive integer $ n_{0} 
< k $ such
 that $ W_{0} \subset V_{0} $ and $ f^{n_{0}}(W_{0}) \subset V_{1} $. This is possible as $f^{n_{0}}$ is a
continuous
function. Next
choose a neighbourhood $ W_{1} $ of $ x_{0} $ and a positive integer $ n_{1} < k $ such
that $ W_{1} \subset W_{0} $ and $
f^{n_{1}}(W_{1}) \subset V_{2} $ . Continuing in this way We have $ k $ such
neighbourhoods of $x_{0}$, namely $ W_{0} \supset
W_{1}  \supset W_{2} \supset W_{3}  \ldots \supset W_{k-1}$. Let $ W$ be their intersection, i.e., the
smallest of them. Then
$ W$ is a neighbourhood of $x_{0} $ such that $ f^{n_{i}}(W) \subset V_{i+1} $ for all $ i=0,1,2, \ldots , k.$
But the $V_{i}^{'}s$
are disjoint. Therefore $ n_{i}s $ have to be distinct. A set of k positive integers each of which is
strictly less than k. This is a
contradiction.

Remark: Under the assumption of the above proposition, f has to be a cyclic permutation
on the finite set. We can restate the above proposition as follows: A dynamical system
$ (X, f ) $
is an $n$-cycle if and only if \\

$ {\\}^{sup}_{V,W} ~~ {\\}^{inf}_{ \ n \ } ~~ \{ n \in N : f^{n}(V) \cap  W $ is  nonempty $ \}~~  < ~~
\infty $

where the supremum is taken over all pairs  $ V, W $ nonempty open sets in $X$. Here the convention is that
the infimum of the empty subset of  $N$ is $ \infty $

{\bf Corollary:} The finite cycles are the only transitive maps satisfying dense periodic points
but not dense large periodic points.

\subsection{Additional remarks}

	The second proof has a distinct advantage over the first. Its method proves some other related
results stated below:

Result-1: Let $ (X,f) $ be an infinite topologically transitive dynamical
system. Let $ \{ F_{n} \} $ be
a sequence of  closed invariant sets such that  

$ { \bigcup }^{ \infty } _{ n = k } F_{n} $ is also dense, for every
positive integer $K$.

Theorem-1 becomes a particular case of Result-1.

Result-2: Let $ (X,f) $ be an infinite topologically transitive dynamical system 
such that the the set of eventually periodic points is dense. Then for every positive
integer $n$, the set of eventually periodic points having an orbit of atleast $n$
elements, is also dense.

Result-2 is a particular case of result-1

Result-2': Let $ (X,F) $ be an infinite topologically transitive dynamical system.
Let $n$ be a positive integer. If each nonempty open set contains an element with
finite orbit, then each nonempty open set contains an element with finite orbit
of size $ \geq  n $

Result-3: Let $ (X, f) $ be topologically transitive and let P be dense. Then
either $f^n$ equals identity for some $n$, or $ P \setminus P_{n}$ is dense for
all $ n \in N $

Result-4: Let  $ (X,f) $ be any dynamical system. Let the set $P$ of periodic points
be a dense set with empty interior. Then $P \setminus P_{n}$ is also dense for
every $n$ in $N$. Note that we are    not  assuming $F$ to be topologically
transitive.

When $ X= R $ or an interval in $R$, we can say more:

Result-5: For every transitive map on $R$ or an interval in $R$, the set $P$ of periodic points
is a dense set with empty interior. Therefore, by Result-4, $ P \setminus P_{n}$
is dense for every $n$.

Result-6: For continuous $ f: I \rightarrow I$ where $P$ is  dense, the only two
possibilities
are

\begin{enumerate}
\item

$f$ is identity or topologically conjugate to a reflection map on some nontrivial
subinterval  $ J \subset I $\\
or
\item

$P \setminus P_{n} $ is dense for every $n$ in $N$.\\
This can be proved by applying a theorem of [ 7 ] that states that $f$ is built from
transitive maps on subinterval. Note that when $f$ is transitive $ (1) $ cannot hold. 
Thus, for all transitive maps on intervals, $(2)$ should hold; of course, this can be
proved
from result-5 also.

\end{enumerate}

{\center {\bf References } }

1.  Devaney, R.L. An introduction to chaotic dynamical systems, Addison Wesley, (1989)

2. Banks, J. Brooks, J. Cairns, G. Davis, G. and  Stacey, P. On Devaney's definition of chaos,
Amer. Math. Monthly, 99, 332-334 (1992). 

3. Block, L.S. and Coppel, W.A. Dynamics in one dimension, Lecture Notes in Mathematics, 1513, Springer, (1992)

4. A. Crannell and M.Martelli, Periodic orbits from Nonperiodic orbits on an interval, Appl. Math. Lett. 
Vol.10, No. 6, 45-47 (1997)

5. Silverman, S. On maps with dense orbits and the definition of chaos, Rocky Mountain Jour. Math. 22,353-375
(1992)

6. Touhey, P. Yet another definition of chaos, Amer. Math. Monthly, 104, 411-414 (1997)

7. Barge, M and Martin, J. Dense periodicity on the interval, Proc. Amer. math. Soc., 94, 731-735 (1985)

\end{document}